\documentclass{elsart_mod}

\usepackage[dvips]{graphicx}
\usepackage{amsmath}
\usepackage{amssymb}
\usepackage{enumerate}
\usepackage{url}

\newcommand{\eref}[1]{(\ref{#1})}

\begin{document}
\begin{frontmatter}
\title{Efficient initial designs for binary response data}
\author{Juha Karvanen}
\address{
Department of Health Promotion and Chronic Disease Prevention,\\ 
National Public Health Institute,\\
Mannerheimintie 166, 00300 Helsinki, Finland\\
juha.karvanen@ktl.fi}

\begin{abstract}
In this paper we introduce a binary search algorithm that efficiently finds initial maximum likelihood estimates for sequential experiments where a binary response is modeled by a continuous factor. The problem is motivated by switching measurements on superconducting Josephson junctions. In this quantum mechanical experiment, the current is the factor controlled by the experimenter and a binary response indicating the presence or the absence of a voltage response is measured. The prior knowledge on the model parameters is typically poor, which may cause the common approaches of initial estimation to fail. The binary search algorithm is designed to work reliably even when the prior information is very poor. The properties of the algorithm are studied in simulations and an advantage over the initial estimation with equally spaced factor levels is demonstrated. We also study the cost-efficiency of the binary search algorithm and find the approximately optimal number of measurements per stage when there is a cost related to the number of stages in the experiment.

\end{abstract}
\end{frontmatter}

\section{Introduction}
Consider an experiment where the probability of success (or failure) is a monotonic function of a factor controlled by the experimenter. The functional form of the response curve is known but the parameters defining its location and slope are unknown. Our goal is to estimate these unknown parameters in the situation where the prior knowledge on the parameters is very poor. Intuitively it is clear that the factor levels should be chosen such a way that both successes and failures are measured as responses -- otherwise the experiment provides very little information on the parameters of the interest. However, if the location and slope of the response curve are unknown, the optimal factor levels are also unknown.

The problem described above is encountered in quantum mechanical experiments called switching measurements \cite{Josephson,single,pekola}. In the experiment, a Josephson junction circuit is placed in a dilution refrigerator at temperature below 4.2~kelvin. The height of a current pulse measured in nanoamperes is the factor controlled  by the experimenter and the absence or presence of a voltage pulse is measured as a response. For a certain narrow range of current pulses, the presence of the voltage response is a random variable depending on the height of the current pulse, which allows the Josephson junctions to be used as ultra-sensitive sensors. Only imprecise limits of the sensitive range can be given in advance: the current must be above zero and a rough upper limit is obtained by measuring the resistance of the circuit in room temperature.

In the literature on the optimal design for binary response models, it is often assumed that prior estimates are readily available or that they can be easily achieved from a pilot experiment. A natural but not always very efficient approach to initial estimation is to choose the initial factor levels that are equally spaced on some initial range. This approach has been employed in several theoretical and practical works, e.g. \cite{Sitter:robustdesign,Sitter:twostagequantal,Chaudhuri:nonlinearexperiments,Chaloner:Bayesianexperimental}. A more advanced approach based on the minimax principle was proposed and studied by Sitter \cite{Sitter:robustdesign} and King and Wong \cite{King:minimaxdoptimal}

In this paper we study the problem of initial design in detail and introduce a binary search algorithm that reliably finds initial maximum likelihood estimates (MLEs) even in the case where the prior information on the parameters is very poor. The binary search algorithm was already successfully applied to switching measurements \cite{optdesign}. We also study the cost-efficiency of the binary search algorithm and find the approximately optimal number of measurements per stage assuming that there is a cost related to the number of stages in the experiment.

The rest of the paper is organized as follows. In Section~\ref{sec:responsemodels} we review some commonly used binary response models and present the related theory of optimal design. The problem of the non-existence of MLEs with small and moderate sample sizes is studied. In Section~\ref{sec:binarysearch}, the binary search algorithm is introduced and its properties are studied. In Section~\ref{sec:equallyspaced}, the binary search algorithm is compared to the initial estimation with equally spaced factor levels. In Section~\ref{sec:costefficient}, we consider a cost model where there is a cost related to the number of measurements and a cost related to the number of stages and find the approximately optimal number of measurements in the binary search. In Section~\ref{sec:application}, we apply the results to switching measurements and discuss the use of binary search in applications similar to the sport fishing example \cite{Sitter:robustdesign,King:minimaxdoptimal}. Section~\ref{sec:conclusion} concludes the paper.

\section{Binary response models} \label{sec:responsemodels}
There are three parametric models, the logit, the probit and the complementary log-log (cloglog) model, that are frequently used to model the dependence between a binary response $Y$ and a continuous factor $x$. All these models can be presented in the framework of generalized linear models \cite{McCullagh:glm}
\begin{equation} \label{eq:binary}
P(Y=1)=\textrm{E}(Y)=F(ax+b),
\end{equation}
where the response curve $F$ is a cumulative distribution function (cdf) and the $a$ and $b$ are the parameters of the model to be estimated. The three response curves commonly used are:
logistic distribution for the logit model
\begin{equation} \label{eq:logit}
F(ax+b)=\frac{\exp(ax+b)}{1+\exp(ax+b)},
\end{equation}
normal distribution for the probit model
\begin{equation} \label{eq:probit}
F(ax+b)=\Phi(ax+b)=\int_{-\infty}^{\infty}\frac{1}{\sqrt{2 \pi}}\exp\left(\frac{-(ax+b)^{2}}{2}\right) \textrm{d} x,
\end{equation}
and the Gompertz distribution for the cloglog model
\begin{equation} \label{eq:cloglog}
F(ax+b)=1-\exp(-\exp(ax+b)).
\end{equation}

The theory of optimal design provides results describing how the factor levels should be chosen in order to estimate parameters $a$ and $b$ optimally. Table~\ref{tab:Doptims} adapted from \cite{Ford:canonicaloptimal} presents the (locally) D-optimal factor levels for the logit, the probit and the cloglog model. The apparent problem discussed in more detail by Minkin \cite{Minkin:optimal} and Sitter \cite{Sitter:robustdesign} is that the optimal factor levels are functions of the unknown parameters that we should estimate. It follows that we can design the optimal experiment only if we already know the parameters we want to estimate, which is an impractical requirement. The problem is often evaded assuming that we already have good initial estimates but less often it is discussed how the good initial estimates are obtained and what follows if initial information is very poor.

\begin{table}[htb]
\caption{D-optimal factor levels for the models~\eref{eq:logit}, \eref{eq:probit} and \eref{eq:cloglog}. The D-optimal factor levels $x_{1}^{*}$ and $x_{2}^{*}$ are solved from equations $ax_{1}^{*}+b=z_{1}^{*}$ and $ax_{2}^{*}+b=z_{2}^{*}$. Columns $F(z_{1}^{*})$ and $F(z_{2}^{*})$ report the cdf levels related to the optimal factors, i.e., the probabilities of response~1. \label{tab:Doptims} }
\begin{center}
\begin{tabular}{rcccc}
model & \multicolumn{4}{c}{D-optimal factors}\\
 & $z_{1}^{*}$ & $z_{2}^{*}$ & $F(z_{1}^{*})$ & $F(z_{2}^{*})$\\ \hline
logit  & -1.543 & 1.543 & 0.176 & 0.824\\
probit  & -1.138 & 1.138 & 0.128 & 0.872\\
cloglog  & -1.338 & 0.980 & 0.231 & 0.930\\
\end{tabular}
\end{center}
\end{table}

Next we study the problem of the non-existence of MLEs. We consider maximum likelihood estimation in small experiments where the sample size $n \leq 200$ and assume that the true values of $a$ and $b$ are known so that the actual (locally) D-optimal design can be used. MLEs cannot be estimated from the measured data if it is possible to classify the responses to 1's and 0's on the basis of the factor levels, e.g. response~1 is obtained if $x>100$ and response~0 is obtained if $x<100$. More formally, the existence of MLEs for the model~\eref{eq:binary} requires that \cite{Albert:mlexistence}
\begin{align} \label{eq:mlexistence}
\max(x_{i}|y_{i}&=0)>\min(x_{i}|y_{i}=1) \textrm{ and } \nonumber \\
\max(x_{i}|y_{i}&=1)>\min(x_{i}|y_{i}=0),
\end{align}
where $y_{i}$ is the measured response for the factor level $x_{i}$.
This condition is fulfilled, for instance, if there are two such factor levels that for each of them both 0's and 1's are measured as responses. It follows that the probability for the non-existence of MLEs under the D-optimal design is equivalent to the probability that the only  0's or 1's are measured as responses in either at $z_{1}^{*}$ or at $z_{2}^{*}$. According to the elementary laws of probability calculus, this probability can be expressed as
\begin{equation}
P(\textrm{No MLE})=q_{1}^{s}+p_{1}^{s}+q_{2}^{s}+p_{2}^{s}-
q_{1}^{s}q_{2}^{s}-q_{1}^{s}p_{2}^{s}-
p_{1}^{s}q_{2}^{s}-p_{1}^{s}p_{2}^{s},
\end{equation}
where $p_{1}=F(z_{1}^{*})$ and $p_{2}=F(z_{2}^{*})$ are read from Table~\ref{tab:Doptims}, $q_{1}=1-p_{1}$, $q_{2}=1-p_{2}$, $s=n/2$ and $n$ is the sample size, which is assumed to be an even number. The probabilities for the non-existence of MLEs in D-optimal designs as a function of sample size are presented in Table~\ref{tab:MLEnonexistence}. It can be seen that surprisingly large sample sizes are needed in order to avoid the risk that MLEs cannot be calculated from the sample. The conclusion here is that even if we are using the D-optimal design, the non-existence of MLEs is a problem that cannot be ignored when the sample size is small. In practical situations where the D-optimal design is not known the problem may be even more serious.

\begin{table}[htb]
\caption{Probabilities for the non-existence of MLEs when the D-optimal factor levels from Table~\ref{tab:Doptims} are used. \label{tab:MLEnonexistence} }
\begin{center}
\begin{tabular}{rccc}
sample size & logit & probit & cloglog\\ \hline
 4  &  0.91584142  &  0.95047339  &  0.95394709 \\
 10  &  0.61553275  &  0.75551315  &  0.77862340 \\
 20  &  0.26765833  &  0.44578006  &  0.52280966 \\
 40  &  0.04117204  &  0.12633768  &  0.23974289 \\
 100  &  0.00012482  &  0.00217817  &  0.02698092 \\
 200  &  0.00000001  &  0.00000237  &  0.00072786\\
 \end{tabular}
\end{center}
\end{table}

\section{Initial estimation with binary search} \label{sec:binarysearch}
We start our introduction to the binary search algorithm considering types of the prior information. There is virtually always some prior information on the parameters available although this prior information may be very inaccurate and difficult to express in mathematical form.
For instance, the following types of prior information may be encountered:
\begin{itemize}
\item Point estimates of the parameters from the previous study/studies.
\item Range of the possible parameter values.
\item Prior distribution of the parameters.
\item Range of sensible factor levels.
\end{itemize}
We concentrate to prior information given as a range of factor levels because it is the simplest type of prior information and the other types can be transformed to this type if needed. It is actually difficult to find a practical situation where no initial interval could be specified. Value zero is often a theoretical lower bound  and some extreme high value can serve as the upper bound. In the other words, we can always specify the initial interval in such a way that it contains all sensible factor levels with high probability. If the prior information is very inaccurate, the initial interval can become very wide.

\begin{figure}
  \includegraphics[width=\columnwidth]{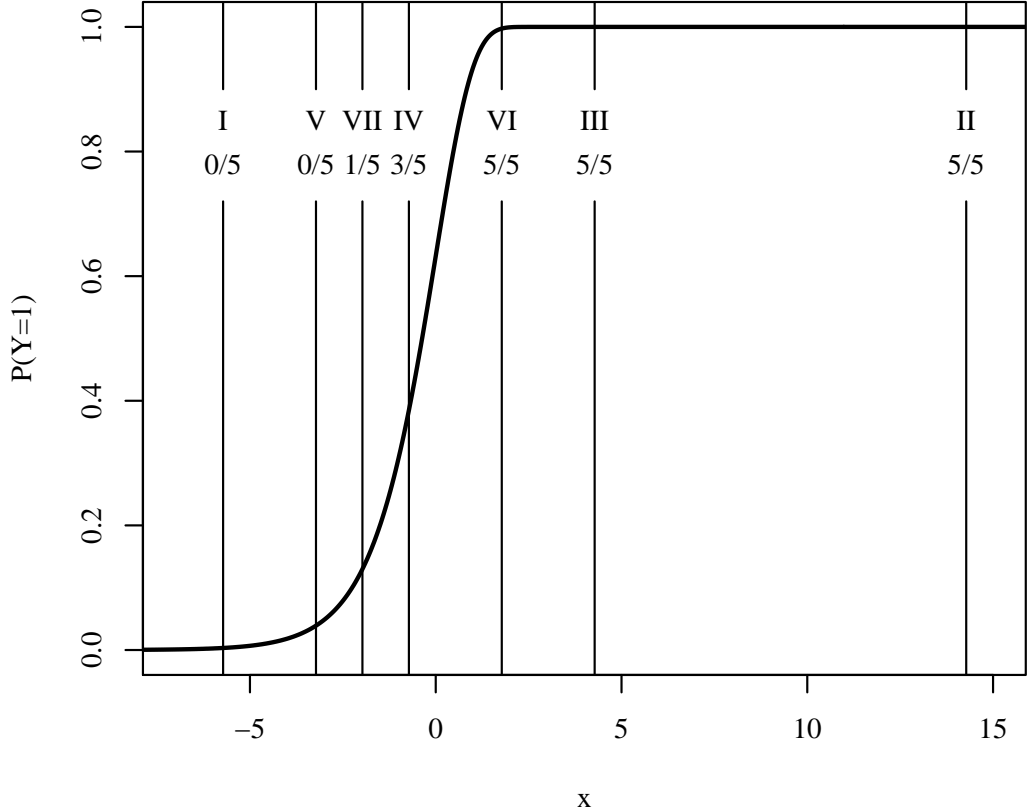}
  \caption{Binary search algorithm in action. The response curve of the cloglog model with $a=1$ and $b=0$ is plotted. The initial interval is defined as $[-5.7,14.3]$.  The vertical lines present the measured factor levels and the Roman numerals indicate the order of the measurements. The numbers below the Roman numerals report the number of responses 1 per the total measurements; e.g. $3/5$ means that three responses 1 were measured out of five measurements. The two points needed for the maximum likelihood estimation were found at measurements IV and VII. The obtained initial estimates of the parameters of the cloglog model and their standard errors are $\hat{a}=1.83\;(0.88)$ and  $\hat{b}=1.88\;(1.29)$.} \label{fig:binarysearch}
\end{figure}

The key idea of the binary search algorithm is to measure at such factor levels that the condition for the existence of MLEs~\eref{eq:mlexistence} is fulfilled as quickly as possible. The motivation for this is that after the MLEs have been found, we can continue the experiment using the design that is optimal for the current estimates of the parameters. The proposed method is similar to the binary search algorithm that finds a root of a continuous function. Instead of a root we are looking for a factor level where both 0's and 1's can be measured as responses. At every step, we measure the responses at the middle point of the current search interval: If only 0's are measured, the middle point is taken as the new starting point of the interval. If only 1's are measured, the middle point is taken as the new end point of the interval. The first part of the binary search ends when both 0's and 1's are measured at same point $x$. Another factor level needed to guarantee the existence of MLEs is then found in the neighborhood of $x$.
The proposed algorithm can be presented as follows:
\begin{enumerate}[1.]
\item Use the previous knowledge to construct such an interval
$[x_{\textrm{min}},x_{\textrm{max}}]$ that we can be sure that
$F(x_{\textrm{max}};a,b)-F(x_{\textrm{min}};a,b)$ is close to 1.
Constructing this kind of interval is usually possible even if we have
very little knowledge on the parameters.
\item Use binary search to find such a point $x$ from the interval
$[x_{\textrm{min}},x_{\textrm{max}}]$ that both 0's and 1's are measured as responses. In binary search, we measure the responses at the middle point of the current interval. If only 0's are measured, the middle point is taken as the new starting point of the interval. If only 1's are measured, the middle point is taken as the new end point of the interval. Let $x$ be found as the
middle point of the interval $[x_{l},x_{u}]$, i.e. $x=(x_{l}+x_{u})/2$.
\item Define $|\epsilon|=(x_{u}-x_{l})/4$. The sign of $\epsilon$ is determined according to the measured response: $\textrm{sign}(\epsilon)=\textrm{sign}(0.5-\bar{y})$, where $\bar{y}$ is the average response for $x$. If $0.5-\bar{y}=0$ the sign of $\epsilon$ is chosen randomly.
\item Measure first at point $x+\epsilon$. If both 0's and 1's are not measured as responses, measure also at point $x-\epsilon$.
\item If MLEs exist for the data measured so far, proceed to the maximum
likelihood estimation. Otherwise, divide $\epsilon$ by two and return to Step 4.
\end{enumerate}
The progress of the algorithm is illustrated in Figure~\ref{fig:binarysearch}. Step 2 requires that we measure at least two times at each point.  When the binary search in Step 2 converges we have found one of the two points required for the maximum likelihood estimation. The other point that is needed is then found in the neighborhood of point $x$ using again binary search.  The procedure works reliably and efficiently because at every step of binary search, the length of the current interval is divided by two. Consequently, the procedure converges exponentially fast regardless of the choice of the initial interval.

After it has been found that condition \eref{eq:mlexistence} holds and the MLEs therefore exist, there are two alternative ways to use the data. We may calculate the MLEs from all data that has been cumulated in the binary search algorithm (MLE method I) or use only the endpoints of the initial interval and the two points found in Steps 2 and 4 of the algorithm (MLE method II). Neither of these approaches is fully satisfactory because the data originate from a procedure with data dependent stopping rules. However, the purpose of the initial design is only to produce some initial estimates for the primary experiment, which means that any sensible estimate based on real data is better than an arbitrary guess.

The performance of the proposed binary search algorithm is studied in simulations where the quality of the prior information varies from poor to very poor. We fix the parameters $a=1$ and $b=0$ and specify the initial interval for factor levels as $[h-d,h]$, where $d>2h$ is the length of the initial interval (a simulation parameter) and $h$ is generated from \mbox{$\textrm{Uniform}(4,d-4)$} distribution. The results can be generalized because arbitrary parameters $a$ and $b$ can be returned to the canonical case $a=1$ and $b=0$ by a linear transformation. Note that the initial intervals are specified in such way that the middle point of the response curve always belongs to the interval but its relative location with respect to the endpoints of the interval is unknown and uniformly distributed. The interval lengths $d \in \{10,15,20,50,100,200,1000\}$ are used in the simulations. The number of measurements per stage $n_{k}$ is also a simulation parameter and has 23 values $\{2,3,4,\ldots,200,500,1000\}$. 5000 simulation runs are performed for each combination of simulation parameters. For each simulation run we record
\begin{itemize}
\item $K$, the number of stages needed for the convergence of the binary search,
\item $\mathbf{J}_{K}$, the Fisher information matrix after the initial design,
\item $\hat{a}$ and $\hat{b}$, the parameter estimates after the initial design.
\end{itemize}
The validity of the initial interval is always first checked measuring the responses at the ends of the initial interval. These measurements are counted as the two first stages and the actual binary search starts from the stage three.

Tables~\ref{tab:bs_cloglogreqn}, \ref{tab:bs_cloglogI} and \ref{tab:bs_cloglogII} report the simulation results for the cloglog model. It can be seen from Table~\ref{tab:bs_cloglogreqn} that more stages are needed when the initial interval widened but the number of the required additional stages is small. In fact, if the length of the initial interval is multiplied by two, one additional stage is required because at each stage of the binary search algorithm, the current interval is halved. It follows that if the length of the initial interval is multiplied by ten, $\log_{2}(10) \approx 3.32$ additional stages are required, which is in agreement with the differences between the observed $K$ for $d=10$, $d=100$, and $d=1000$ in Table~\ref{tab:bs_cloglogreqn}.

The estimates using the MLE method I and the MLE method II are presented in Tables~\ref{tab:bs_cloglogI} and \ref{tab:bs_cloglogII}, respectively.
The averages of estimates $\hat{a}$ and $\hat{b}$ indicate that $a$ and $b$ are systemically overestimated with the method I.  The bias is large for small stage sizes (average $\hat{a}=2.56$ for $n_{k}=3$) and small for large stage sizes (average $\hat{a}=1.04$ for $n_{k}=100$) but it seems to be independent of the length of the initial interval $d$.  For $n_{k}=2$, the average bias is very large indicating that some estimates were very large.  The reason for the bias is that there is only limited number of outcomes when the initial estimation may stop. With method II, $a$ is underestimated when $n_{k}=2$.   On the basis of these results, it is clear that the stage size $n_{k}=2$ should not be used. With method II, parameter $a$ is slightly overestimated for $n_{k}>2$ when $d=10$ and underestimated when $d=100$ and $d=1000$. Parameter $b$ is slightly underestimated. Similarly to method I, the bias decreases when $n_{k}$ increases.

There are several ways to deal with the bias. The safest way is to choose large $n_{k}$. If small $n_{k}$ is preferred, the results in Tables~\ref{tab:bs_cloglogI} and \ref{tab:bs_cloglogII} could be employed in making $\hat{a}$ and $\hat{b}$ unbiased but this is not necessarily a good solution for method I because the distribution of $\hat{a}$ is strongly skewed, which can be seen comparing the mean and the 95th percentile. Simpler approach is to use the estimates without a bias correction and expect that the sequential design will quickly correct the bias. An overestimated $a$ causes the D-optimal factor levels calculated according to Table~\ref{tab:Doptims} shift towards the center of the response curve. This is not particularly harmful because both 0's and 1's can be measured as responses in the central area of the response curve.

The Fisher information after convergence of the binary search is almost independent of the length of the initial interval. This can be explained by the fact that when the initial interval is wide, the measurement points at the first stages fall to the tails of the response curve and thus provide very little information in terms of the Fisher information. The number of stages that provide the information is therefore nearly constant. The same fact explains why the Fisher information is almost the same for method I and method II.

\begin{table}[htb]
\caption{Simulation results on the number of stages needed for the convergence using the binary search algorithm under the cloglog model.
The numbers reported are the mean, the 5th percentile and the 95th percentile calculated from 5000 simulation runs for different lengths of the initial interval $d$ and number of measurements per stage $n_{k}$.\label{tab:bs_cloglogreqn} }
\begin{center}
\begin{small}
\begin{tabular}{rrrrrrrrrrr}
& \multicolumn{3}{c}{$d=10$} & \multicolumn{3}{c}{$d=100$} & \multicolumn{3}{c}{$d=1000$}\\
$n_{k}$ & mean & 5\% & 95\% & mean & 5\% & 95\% & mean & 5\% & 95\%\\
 2      &  7.53  &  5 &  11  &  10.93    &  8  &  14   &  14.26  &  12  &  18 \\
 3      &  6.53  &  4 &  9   &  9.87     &  8  &  13   &  13.27  &  11  &  16 \\
 4      &  6.09  &  4 &  8   &  9.52     &  8  &  12   &  12.93  &  11  &  16 \\
 5      &  5.81  &  4 &  8   &  9.25     &  7  &  12   &  12.69  &  11  &  16 \\
 10     &  5.11  &  4 &  7   &  8.72     &  7  &  11   &  12.19  &  10  &  15 \\
 100    &  4.25  &  4 &  5   &  7.83     &  6  &  10   &  11.35  &  10  &  14 \\
 1000   &  4.09  &  4 &  5   &  7.28     &  6  &  9    &  10.91  &  9   &  14
\end{tabular}
\end{small}
\end{center}
\end{table}

\begin{table}[htb]
\caption{Simulation results for the cloglog model using the binary search algorithm and MLE method I.
The numbers reported are the mean, the 5th percentile and the 95th percentile calculated from 5000 simulation runs for different lengths of the initial interval $d$ and number of measurements per stage $n_{k}$.\label{tab:bs_cloglogI} }
\linespread{0.8}
\selectfont
\begin{center}
\begin{small}
\begin{tabular}{rrrrrrrrrrr}
\\
\multicolumn{10}{l}{$\hat{a}$, estimated slope parameter $(a=1)$}\\
& \multicolumn{3}{c}{$d=10$} & \multicolumn{3}{c}{$d=100$} & \multicolumn{3}{c}{$d=1000$}\\
 $n_{k}$ & mean & 5th & 95th & mean & 5th & 95th & mean & 5th & 95th\\
 2  &  199.58  &  0.50  &  28.79  &  221.84  &  0.64  &  27.20  &  23.22  &  0.58  &  33.49 \\
 3  &  2.56  &  0.60  &  5.11  &  21.03  &  0.50  &  7.17  &  8.30  &  0.64  &  5.86 \\
 4  &  1.69  &  0.57  &  3.98  &  1.64  &  0.55  &  3.64  &  1.65  &  0.58  &  3.56 \\
 5  &  1.47  &  0.63  &  3.00  &  1.46  &  0.57  &  2.85  &  1.46  &  0.57  &  2.95 \\
 10  &  1.20  &  0.67  &  2.07  &  1.19  &  0.68  &  2.03  &  1.19  &  0.61  &  2.05 \\
 100  &  1.04  &  0.85  &  1.29  &  1.02  &  0.78  &  1.27  &  1.01  &  0.79  &  1.25 \\
 1000  &  1.00  &  0.94  &  1.07  &  0.99  &  0.83  &  1.11  &  1.00  &  0.86  &  1.15\\
\\
\multicolumn{10}{l}{$\hat{b}$, estimated location parameter $(b=0)$}\\
 & \multicolumn{3}{c}{$d=10$} & \multicolumn{3}{c}{$d=100$} & \multicolumn{3}{c}{$d=1000$}\\
  $n_{k}$ & mean & 5th & 95th & mean & 5th & 95th & mean & 5th & 95th\\
 3  &  0.67  &  -1.05  &  2.76  &  40.77  &  -1.05  &  2.74  &  14.72  &  -1.01  &  2.69 \\
 4  &  0.27  &  -0.89  &  1.78  &  0.29  &  -0.85  &  1.83  &  0.26  &  -0.83  &  1.76 \\
 5  &  0.20  &  -0.76  &  1.47  &  0.20  &  -0.75  &  1.49  &  0.19  &  -0.76  &  1.46 \\
 10  &  0.09  &  -0.56  &  0.92  &  0.08  &  -0.58  &  0.86  &  0.08  &  -0.57  &  0.86 \\
 100  &  0.02  &  -0.21  &  0.26  &  0.01  &  -0.26  &  0.29  &  -0.02  &  -0.30  &  0.25 \\
 1000  &  0.00  &  -0.07  &  0.07  &  -0.03  &  -0.26  &  0.11  &  0.01  &  -0.15  &  0.22\\
\\
\multicolumn{10}{l}{$D$, the square root of the determinant of the Fisher information}\\
& \multicolumn{3}{c}{$d=10$} & \multicolumn{3}{c}{$d=100$} & \multicolumn{3}{c}{$d=1000$}\\
$n_{k}$ & mean & 5th & 95th & mean & 5th & 95th & mean & 5th & 95th\\
 2  &  1.73  &  0.08  &  4.96  &  1.66  &  0.09  &  3.75  &  1.64  &  0.07  &  4.37 \\
 3  &  2.69  &  0.51  &  6.44  &  2.66  &  0.51  &  6.79  &  2.73  &  0.55  &  6.09 \\
 4  &  3.76  &  1.08  &  7.60  &  3.61  &  1.10  &  7.60  &  3.67  &  0.94  &  8.24 \\
 5  &  4.56  &  1.30  &  9.45  &  4.46  &  1.52  &  9.13  &  4.53  &  1.51  &  9.52 \\
 10  &  8.33  &  3.25  &  15.21  &  8.35  &  3.95  &  14.76  &  8.34  &  3.80  &  15.03 \\
 100  &  68.58  &  30.98  &  101.87  &  55.56  &  20.57  &  89.39  &  55.60  &  24.47  &  89.33 \\
 1000  &  679.14  &  425.46  &  952.25  &  370.71  &  103.53  &  756.21  &  323.80  &  86.97  &  658.73
\end{tabular}
\end{small}
\end{center}
\end{table}

\begin{table}[htb]
\caption{Simulation results for the cloglog model using the binary search algorithm and MLE method II.
The numbers reported are the mean, the 5th percentile and the 95th percentile calculated from 5000 simulation runs for different lengths of the initial interval $d$ and number of measurements per stage $n_{k}$.\label{tab:bs_cloglogII} }
\linespread{0.8}
\selectfont
\begin{center}
\begin{small}
\begin{tabular}{rrrrrrrrrrr}
\\
\multicolumn{10}{l}{$\hat{a}$, estimated slope parameter $(a=1)$}\\
& \multicolumn{3}{c}{$d=10$} & \multicolumn{3}{c}{$d=100$} & \multicolumn{3}{c}{$d=1000$}\\
 $n_{k}$ & mean & 5th & 95th & mean & 5th & 95th & mean & 5th & 95th\\
 2  &  0.86  &  0.40  &  1.70  &  0.24  &  0.09  &  0.60  &  0.06  &  0.01  &  0.18 \\
 3  &  1.21  &  0.44  &  3.19  &  0.74  &  0.10  &  2.55  &  0.78  &  0.01  &  2.04 \\
 4  &  1.04  &  0.44  &  2.52  &  0.79  &  0.11  &  2.01  &  0.73  &  0.02  &  1.80 \\
 5  &  1.02  &  0.47  &  2.26  &  0.86  &  0.11  &  2.53  &  0.80  &  0.02  &  2.02 \\
 10  &  1.03  &  0.53  &  1.87  &  0.91  &  0.23  &  1.74  &  0.90  &  0.19  &  1.73 \\
 100  &  1.02  &  0.82  &  1.28  &  0.98  &  0.71  &  1.26  &  0.96  &  0.67  &  1.24 \\
 1000  &  1.00  &  0.94  &  1.07  &  0.98  &  0.81  &  1.11  &  0.99  &  0.84  &  1.14\\
\\
\multicolumn{10}{l}{$\hat{b}$, estimated location parameter $(b=0)$}\\
 & \multicolumn{3}{c}{$d=10$} & \multicolumn{3}{c}{$d=100$} & \multicolumn{3}{c}{$d=1000$}\\
  $n_{k}$ & mean & 5th & 95th & mean & 5th & 95th & mean & 5th & 95th\\
 2  &  -0.02  &  -0.96  &  1.06  &  -0.30  &  -1.17  &  0.36  &  -0.39  &  -1.23  &  0.33 \\
 3  &  0.09  &  -1.01  &  1.20  &  -0.13  &  -1.47  &  1.03  &  -0.11  &  -1.65  &  1.06 \\
 4  &  0.01  &  -0.92  &  1.12  &  -0.11  &  -1.16  &  1.02  &  -0.15  &  -1.23  &  0.98 \\
 5  &  -0.00  &  -0.91  &  1.03  &  -0.07  &  -1.26  &  1.01  &  -0.12  &  -1.45  &  1.01 \\
 10  &  0.00  &  -0.73  &  0.80  &  -0.08  &  -1.00  &  0.70  &  -0.09  &  -1.02  &  0.71 \\
 100  &  -0.00  &  -0.26  &  0.26  &  -0.02  &  -0.37  &  0.27  &  -0.07  &  -0.54  &  0.24 \\
 1000  &  -0.00  &  -0.07  &  0.07  &  -0.05  &  -0.30  &  0.10  &  -0.00  &  -0.22  &  0.22\\
\\
\multicolumn{10}{l}{$D$, the square root of the determinant of the Fisher information}\\
& \multicolumn{3}{c}{$d=10$} & \multicolumn{3}{c}{$d=100$} & \multicolumn{3}{c}{$d=1000$}\\
$n_{k}$ & mean & 5th & 95th & mean & 5th & 95th & mean & 5th & 95th\\
 2  &  1.98  &  0.22  &  4.71  &  3.28  &  0.28  &  9.77  &  6.75  &  0.34  &  27.02 \\
 3  &  2.84  &  0.42  &  6.70  &  3.87  &  0.52  &  12.32  &  6.94  &  0.65  &  32.24 \\
 4  &  3.76  &  0.91  &  8.70  &  4.20  &  0.87  &  12.92  &  6.80  &  0.78  &  30.65 \\
 5  &  4.42  &  1.08  &  9.84  &  4.42  &  1.17  &  13.16  &  6.59  &  1.14  &  28.09 \\
 10  &  7.72  &  2.96  &  15.22  &  6.81  &  2.75  &  13.10  &  7.14  &  2.76  &  12.33 \\
 100  &  65.10  &  27.30  &  101.56  &  50.68  &  17.02  &  83.85  &  50.03  &  17.46  &  80.59 \\
 1000  &  670.40  &  396.07  &  952.25  &  361.69  &  87.16  &  751.10  &  314.30  &  77.38  &  652.07
\end{tabular}
\end{small}
\end{center}
\end{table}

Full simulation results including results for the logit model and the probit model can be obtained from the author by a request. The results for the logit model and the probit model are otherwise similar to the results for the cloglog model but the estimates of $b$ are unbiased because of the symmetry of the response curves.

\section{Comparison with equally spaced factors} \label{sec:equallyspaced}
A commonly used practice for initial design is to choose factor levels that are equally spaced in the design space. Typically, the number of levels is selected on ad-hoc basis. More systematic approaches were proposed by Sitter \cite{Sitter:robustdesign} and King and Wong \cite{King:minimaxdoptimal} who selected the initial design using the minimax principle. It was found by them that the number of the initial levels should be increased when the uncertainty about the model parameters increases. However, when the initial range is wide and the number of initial levels is large, most of the levels will lie far from the D-optimal levels and thus provide very little information on the parameters of the interest. When the prior information is very poor, there is a substantial risk that equal spacing will just utilize resources without any guarantee of finding the MLEs. We study the performance of initial estimation with equally spaced levels in the simulation setting presented in Section~\ref{sec:binarysearch}. The additional simulation parameters are the number of levels $M$ and the number of measurements per each level $n_{m}$. $M$ has values 5, 10, 50, 100 and 1000 and $n_{m}$ has some suitably chosen values that lead to a reasonable sample size $M n_{m}$. For each combination of $d$, $M$ and $n_{m}$, 10000 simulation runs are performed  and the existence of MLEs is checked using the condition~\eref{eq:mlexistence}.

Table~\ref{tab:sp_cloglog} reports the empirical probabilities for the non-existence of MLEs under the cloglog model. The results show a dramatic change in the performance when the quality of the prior information decreases. For instance, if $d=10$, the strategy with $M=10$ and $n_{m}=10$ finds MLEs practically always but if $d=100$ the same strategy fails almost without exceptions. If $d=1000$, the existence of MLEs cannot be guaranteed even if $M=1000$ and $n_{m}=2$. In general it seems that if the sample size $M n_{m}$ is fixed, it is a better to have large $M$ and small $n_{m}$ than vice versa. This is in agreement with the results in \cite{Sitter:robustdesign} and \cite{King:minimaxdoptimal}.

Initial estimation with equally spaced levels cannot compete with the binary search approach. Comparing Tables~\ref{tab:bs_cloglogreqn} and \ref{tab:sp_cloglog} we see   that the required total sample sizes are significantly smaller in the binary search. By looking the row corresponding to $n_{k}=5$ in Table~\ref{tab:bs_cloglogreqn} we notice that using the binary search the sample size of 80 ($=5 \cdot 16$) was sufficient for finding the MLEs even when $d=1000$. In Table~\ref{tab:sp_cloglog}, not even the sample size of 2000 ($=1000 \cdot 2$) was sufficient to eliminate the risk of MLE non-existence.

\begin{table}[htb]
\caption{Empirical probabilities of MLE non-existence for the cloglog model in initial estimation with equally spaced factor levels. The simulation parameters are $d$, the length of the initial interval, $M$, the number of factor levels, and $n_{m}$, the number of measurements per each level. The results are means from 10000 simulation runs. \label{tab:sp_cloglog} }
\begin{center}
\begin{tabular}{rrccc}
  $M$  &  $n_{m}$  &  $d=10$  &  $d=100$  &  $d=1000$ \\ \hline
 5  &  1  &  0.9686  &  1.0000  &  1.0000 \\
 5  &  2  &  0.9084  &  1.0000  &  1.0000 \\
 5  &  10  &  0.4641  &  1.0000  &  1.0000 \\
 5  &  20  &  0.2165  &  1.0000  &  1.0000 \\
 5  &  100  &  0.0029  &  1.0000  &  1.0000 \\
 5  &  1000  &  0.0000  &  1.0000  &  1.0000 \\
 10  &  1  &  0.7609  &  1.0000  &  1.0000 \\
 10  &  2  &  0.4746  &  1.0000  &  1.0000 \\
 10  &  5  &  0.1000  &  1.0000  &  1.0000 \\
 10  &  10  &  0.0058  &  0.9998  &  1.0000 \\
 10  &  50  &  0.0000  &  1.0000  &  1.0000 \\
 10  &  500  &  0.0000  &  0.9955  &  1.0000 \\
 50  &  1  &  0.0164  &  0.9407  &  1.0000 \\
 50  &  2  &  0.0000  &  0.8354  &  1.0000 \\
 50  &  20  &  0.0000  &  0.0756  &  1.0000 \\
 100  &  1  &  0.0000  &  0.7148  &  1.0000 \\
 100  &  2  &  0.0000  &  0.4014  &  1.0000 \\
 100  &  10  &  0.0000  &  0.0038  &  0.9997 \\
 1000  &  1  &  0.0000  &  0.0003  &  0.7070 \\
 1000  &  2  &  0.0000  &  0.0000  &  0.4024
\end{tabular}
\end{center}
\end{table}

\section{Cost-efficient initial designs} \label{sec:costefficient}
 In this section we calculate cost-efficient stage sizes for the proposed binary search procedure. It is assumed that there is a cost related to the number of stages and a cost related to the number of observations. Without loss of generality, we fix the marginal cost of making one additional measurement to one and denote the cost of having one additional stage as $C_{S}$. The same cost scheme was used for D-optimal sequential design in \cite{costdesign}.
 Our goal is to find the optimal number of measurements for initial estimation with the binary search algorithm as a function of the stage cost $C_{S}$ and the length of the initial interval $d$. Optimality is defined here as the cost needed for finding the initial MLEs. We restrict only to initial designs were the number of measurements is the same for each stage of the binary search. If initial estimation with stage size $n_{k}$ takes $K$ stages, the total cost becomes
 \begin{equation} \label{eq:costmodel}
 C(n_{k})=K C_{S}+K n_{k}.
 \end{equation}
 When the stage cost $C_{S}$ is fixed, we can find the stage size $n_{k}$ with the smallest total cost $ C(n_{k})$ among the stage sizes used in the simulation in Section~\ref{sec:binarysearch}. Using this approach the optimal stage size $n_{k}$ was calculated when the stage cost $C_{S}$ took 39 values from 0.001 to 1000. Based on these data, the following model was build for the the optimal $n_{k}$
 \begin{equation} \label{eq:optnk}
\log(n_{k})=\alpha+\beta \log d+\gamma \log C_{S},
 \end{equation}
where the estimated parameters $\alpha$, $\beta$, and $\gamma$ for the logit, the probit and the cloglog model are presented in Table~\ref{tab:models}. In practical experiments, $d$ is always unknown a priori and $C_{S}$ may be known or unknown. We recommend that the choice of $n_{k}$ is based on the optimal $n_{k}$ calculated for several combinations of $d$ and $C_{S}$ using equation~\eref{eq:optnk}.

\begin{table}[htb]
\linespread{1}
\selectfont
\caption{Summary of estimated models~\eref{eq:optnk}.  \label{tab:models} }
\begin{center}
\begin{tabular}{lcc}
\textbf{logit} & \multicolumn{2}{c}{$R^{2}=0.9613$}\\
& estimate & standard error\\
Intercept $\alpha$           &          -0.964  &  0.069 \\
$\log d$ $\beta$             &           -0.200  &  0.013 \\
$\log C_{S}$ $\gamma$        &           0.839  &  0.012 \\
\\
\textbf{probit}  & \multicolumn{2}{c}{$R^{2}=0.9765$}\\
& estimate & standard error\\
Intercept $\alpha$           &       -1.27  &  0.058 \\
$\log d$ $\beta$             &       -0.223  &  0.011 \\
$\log C_{S}$ $\gamma$        &       0.892  &  0.010 \\
\\
\textbf{cloglog} & \multicolumn{2}{c}{$R^{2}=0.9753$}\\
& estimate & standard error\\
Intercept $\alpha$           &       -1.063  &  0.057 \\
$\log d$ $\beta$             &       -0.214  &  0.010 \\
$\log C_{S}$ $\gamma$        &       0.850  &  0.010
\end{tabular}
\end{center}
\end{table}

\section{Applications} \label{sec:application}
We reconsider the measurement data from \cite{optdesign} in order to illustrate the application of the binary search algorithm. In the experiment, a sample consisting of aluminium--aluminium oxide--aluminium Josephson junction  circuit in a dilution refrigerator at 20~millikelvin temperature was connected to computer controlled measurement electronics in order to apply the current pulses and record the resulting voltage pulses. The resistance of the sample at room temperature suggested that a pulse of 300~nA always causes a switching (response~1), which gave the upper limit for the
initial estimation. The lower limit for the initial estimation, 200~nA was roughly estimated from the dimensions of the Josephson junction by an experienced physicist. Alternatively, 0~nA may be used as the lower limit for the initial estimation. The initial estimation used the stage size 50.
A major difference compared to the algorithm presented in Section~\ref{sec:binarysearch} was that the starting value of $\epsilon$ was set as one. In this particular situation, this was a good choice because the MLEs were found after measuring only at four points. Taking the final estimates $\hat{a}=0.240$ and $\hat{b}=-60.628$ calculated from 117288 data points as the true parameters, we study the experiment in the light of the simulation results of Sections~\ref{sec:binarysearch} and \ref{sec:costefficient}. The standardized length of the initial interval is $d=0.240 \cdot (300-200)=24$, or $d=0.240 \cdot 300=72$ for the wider interval. The standardized stage cost $C_{S}$ was estimated from the data and the value $C_{S}=228.4$ was obtained \cite{costdesign}. The model~\eref{eq:optnk} gives the optimal $n_{k}=\exp(-1.063-0.214\log(24)+0.85 \log(228.4)) \approx 18$ for the initial interval $[200,300]$, and $n_{k} \approx 14$ for the initial interval $[0,300]$. According to Table~\ref{tab:bs_cloglogreqn} we expect that the number of stages needed is between 4 and 11.

In addition to switching measurement, the binary search algorithm may be useful also in some other problems. Sitter \cite{Sitter:robustdesign} and King and Wong \cite{King:minimaxdoptimal} consider an application where the problem was to evaluate the economic value of the sport fishing in British Columbia tidal waters. Fishermen were asked the question ``If the cost of your fishing trip had been $x$ dollars higher today, would you still have gone fishing?'' Originally, $x$ varied from \$1 to \$50 but in the middle of the survey it was found that a number of fishermen were willing to pay even more than \$50 extra and the range of $x$ was widened to cover values from \$1 to \$100. Obviously, the original survey was not designed optimally. As an improvement, Sitter \cite{Sitter:robustdesign} proposed a robust seven point design with equal weights for this problem and King and Wong \cite{King:minimaxdoptimal} improved Sitter's results by proposing a robust nine point design with unequal weights. This and similar surveys could benefit from sequential designs with efficient initial design. By using a cellular phone the interviewer could have a direct contact to the computer server collecting the data, which would allow the use of the binary search algorithm described in this paper. In this scenario, the interviewer sends the response of a fisherman to the server that after short calculations returns the amount of dollars for the next question. This would minimize the loss of efficiency due to badly chosen design space.

\section{Conclusion} \label{sec:conclusion}
When the problem of optimal design in binary response experiments is considered in statistical literature it is often assumed that the initial estimates are available. In this paper we propose a binary search algorithm that quickly and reliably returns initial estimates. After finding the initial MLEs the experiment may continue with likelihood based sequential estimation.

Results provided in the paper consider the logit, the probit and the cloglog model but the proposed approach is not restricted to any particular parametric model. On the other hand, the proposed approach relies on the idea of sequential design. If there is a considerable delay between choosing the factor level and measuring the response, sequential designs may be impractical. This applies to many biomedical dose-response trials. In contrast, in many physical experiments the response is measured almost instantly but the total time available for the experiment may be restricted. Efficient sequential designs may be very useful in these problems.

The presented simulations compared the binary search algorithm to the initial estimation with equally spaced factor levels.
The key difference between the approaches is that with equally spaced levels, a good guess for the initial interval must be available whereas with the binary search algorithm, any initial interval that contains the middle point of the response curve is sufficient. If the prior information is very poor, the proposed binary search approach has a huge advantage over the initial estimation with equally spaced levels. The importance of reliable initial estimation is emphasized by the fact that if the prior information is poor it is often difficult to deduce how poor it exactly is.

\section*{Acknowledgements}
The author thanks the Associate Editor and Dr. Juha J. Vartiainen for useful comments.

\bibliographystyle{elsart-num}

\end{document}